# Topology of Vibro-Impact Systems in the Neighborhood of Grazing

Sergey Kryzhevich[1]

*Faculty of Mathematics and Mechanics, Saint-Petersburg State University, Russia,*

Marian Wiercigroch

*Centre for Applied Dynamics Research, The University of Aberdeeen, UK*

**Abstract.** The grazing bifurcation is considered for the Newtonian model of vibro-impact systems. A brief review on the conditions, sufficient for existence of a grazing family of periodic solutions, is given. The properties of these periodic solutions are discussed. A plenty of results on the topological structure of attractors of vibro-impact systems is known. However, since the considered system is strongly nonlinear, these attractors may be invisible or, at least, very sensitive to changes of parameters of the system. On the other hand, they are observed in experiments and numerical simulations. We offer (Theorem 2) an approach which allows to explain this contradiction and give a new robust mathematical model of the non-hyperbolic dynamics in the neighborhood of grazing.

[1]Corresponding author. Tel.: +79219181076.
*E-mail address:* kryzhevitz@rambler.ru, kryzhevich@hotmail.com (S. Kryzhevich).

# 1. Introduction

We study nonlinear systems, describing impact dynamics. One of the most interesting phenomena, occurring in impact systems, is the so-called grazing, where impact velocity tends to zero. The stability of the impacting motion was first comprehensively considered by Peterka [1]. A decade later Thompson and Ghaffari [2] conducted numerical studies showing complexity of dynamical responses for a simple impacting system. Mathematically, a grazing bifurcation corresponds to the case when there exists a family of periodic solutions, continuously depending on a parameter of the system and such that the normal component of the impact velocity vanishes at the bifurcation point.

There is a rich literature on analysis of the dynamics near grazing in a general class of impact oscillator systems, most notably the comprehensive work of Nordmark [3–10]; see [11,12] for an overview. A powerful technique of analysis is to derive a so-called discontinuity mapping which has a square-root singularity, close to the impacting orbit. This can then be combined with an analytic Poincaré map to give a so-called grazing normal form whose dynamics can be shown to be topologically equivalent to those of the underlying flow. See also the works of Whiston [13–15], Budd, Dux [16,17] and Chillingworth [18,19] for some geometrical analysis of the underlying strange attractors in impact oscillators. The structure of the unstable manifold of a limit cycle near grazing was studied in the paper [19].

In this paper we have a different goal, namely to show that the topological features of the resulting attractors can be characterized in terms of the standard tools for chaotic dynamics. In particular, we use the following definition of chaos, given by Devaney [20].

**Definition 1.** Let $U, V \in M$ be domains of a smooth manifold. Consider a diffeomorphism $T : U \to V$. An infinite compact hyperbolic invariant set $K \subset U \bigcap V$ is called *chaotic* if the following conditions are satisfied.

1. The periodic points of the mapping $T$ are dense in $K$.

2. There is a point $p_0 \in K$, whose orbit
$$O(p_0) = \{T^k(p_0) : k \in \mathbb{Z}\}$$
is dense in $K$.

The periodic points of a chaotic invariant set must be unstable. An attractor is called *strange* if it contains a chaotic invariant set.

It was shown by Nordmark [3,6] that the periodic solutions must be unstable in a small neighborhood of grazing. The behavior of the corresponding Lyapunov exponents was studied. It was shown in many papers (experimentally [22–26], numerically [27–30] and analytically [6], [31–34] that the grazing is one of typical reasons for chaotic behavior in vibro-impact systems (VIS).

Usually, the Devaney chaos is found via existence of a transversal homoclinic point, corresponding to a transversal intersection of the stable and the unstable manifold of a fixed point (the Smale-Birkhof theorem, [35]).

An approach to study the structure of the stable and unstable manifold was offered in the article [28] (see also [29]). The homoclinic point is obtained by calculation of a bent of the unstable manifold in the points of non-smoothness. This approach was applied for VIS in the current paper (Section 4).

The other types of chaos are also possible for vibro-impact systems. The sufficient conditions for the Li-Yorke chaos in vibro-impact systems have been given in [11,30] and [36], the stochastic properties of impact oscillators were studied in [37,38].

However, due to the strong nonlinearity of vibro-impact systems, especially the non-s.d.f. ones, the ratios of Lyapunov exponents may be very big for a near-grazing periodic solution. This follows from the normal form of the square root singularity [6,11]. Hence, the corresponding strange attractor may be very sensitive to the changes of parameters of the systems. It may become invisible [25] i.e. it does not physically occur but could occur if the appropriate parameters and initial data were chosen.

Roughly speaking, we have an object (a part of an attractor) which looks like a periodic point on the computer screen. On the other hand the structure of invariant manifolds of this object may be complex, they may look like foliations over the Cantor set. Then, the following questions arise.

1. Does the inner structure of an invisible attractor have an influence on the global behavior of solutions?

2. How can the invariant manifolds, corresponding to subsets of an invisible attractor look like?

This problem has been analyzed by one of the co-authors [25,26]. For a linear harmonic oscillator with an elastic impact the structure of the unstable manifold of the periodic point was studied (Fig. 1). The following system is considered:

$$\dot{x} = y;$$
$$\dot{y} = a\omega^2 \sin(\omega t) - 2\xi y - x - \beta(x-e)H(x-e).$$

The last term at the second equation corresponds to the soft impact (the stiffness of the delimiter (see Fig. 2) is much larger than one of the spring). Both the attractors at the Fig. 1 correspond to the values $\xi = 0.01$ (damping), $e = 1.26$ (position of the delimiter), $a = 0.7$ (excitation amplitude), $\beta = 29$ (stiffness ratio), $H$ is the Heaviside function. The considered grazing bifurcation corresponds to $\omega = 0.801928$.

The left hand side of the figure shows an unstable manifold of a 3-periodic point, $\omega = 0.802$ the right hand side shows one of a chaotic invariant set, $\omega = 0.8023$, which occurs for a very narrow band of the parameter $\omega$ (see [25] for details). One can see that whatever the considered invariant set is, the corresponding unstable manifold has a visible complex structure. So, the structure of an unstable manifold or, more precisely, one of its numerical approximation, traces one of the corresponding invariant set, but is much more easy to be studied numerically.

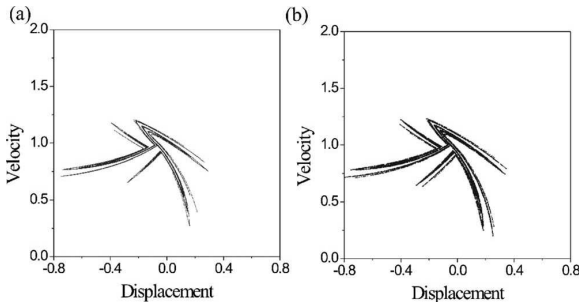

*Fig. 1.*

In the current paper, we give an attempt to explain such phenomena and offer an approach, which gives a new type of chaotic dynamics. Instead of considering the dynamics of points, we introduce (Theorem 2) the dynamical system on an infinite set of disks of codimension 1. We show that this dynamics is transitive and periodic disks are dense. Constructing these so-called "admissible disks" approximating the invariant manifolds, simplifies the studying of the inner structure of the invariant set.

The obtained dynamics is structurally stable. Unlike the classical results of chaotic dynamics (like the Smale-Birkhof theorem on the chaotic dynamics in a neighborhood of a transversal homoclinic point), we do not assume that all the Lyapunov exponents of the bifurcating periodic/fixed point are non-zero in a neighborhood of the bifurcation value. Comparing with the similar results from the theory of partial hyperbolicity [39–47] (see also references therein), we do not need any special assumptions like dynamical coherence (integrability of the central foliation) or dimension 1 of the central manifold.

This dynamics of disks may be considered as a new model of the near-grazing oscillations.

The rest of the paper is organized as follows. At the Section 2 we introduce the Newtonian model of a vibro-impact system and discuss the basic properties of the corresponding solutions. At the Section 3 we define a grazing family of periodic solutions and analyze, when such a family exists. The Section 4 is the main. There, we provide some conditions, sufficient for existence of a non-hyperbolic homoclinic point in a neighborhood of grazing. The existence of such a point does not imply the classical Devaney chaos, so we introduce a special dynamical system on the space of the so-called admissible disks which is described by the symbolical dynamics. A simple example, illustrating the main result of the paper, is considered at the Section 5. At the Section 6 we discuss how the obtained results may be applied. The conclusions are given at the Section 7.

Some definitions and results are accompanied with remarks where the corresponding statements are explained at a less formal level.

## 2. Mathematical model of an impact oscillator.

*2.1. General assumptions and equations of motion of the free flight*

Let $n \in \mathbb{N}$ be the number of degrees of freedom of the considered motion. Consider the half-space

$$\Lambda = [0, +\infty] \times \mathbb{R}^{2n-1}.$$

Denote by $|\cdot|$ the Euclidean norm. Let $\text{col}(a_1, \ldots, a_m)$ be the column vector consisting of the elements $a_1, \ldots, a_m$. For any vector

$$b = \text{col}(b_1, \ldots, b_n)$$

we define $b_{tan} = (b_2, \ldots, b_n)$.

Suppose that the free flight motion between impacts is described by the following equations

$$\ddot{x}_k = f_k(t, x_1, \dot{x}_1, \ldots, x_n, \dot{x}_n, \mu), \qquad k = 1, \ldots, n.$$

or, equivalently, by the system

$$\dot{x}_k = y_k; \qquad \dot{y}_k = f_k(t, z, \mu), \quad k = 1, \ldots, n. \qquad (1)$$

This system may be rewritten in the form

$$\dot{z} = F(t, z, \mu) = \mathrm{col}(y_1, f_1(t, z, \mu), \ldots, y_n, f_n(t, z, \mu)),$$
$$z = \mathrm{col}(x_1, y_1, \ldots, x_n, y_n).$$

We suppose that the parameter $\mu$ is scalar. The continuous function $F : \mathbb{R}^{2n+2} \to \mathbb{R}^{2n}$ is assumed to be $T$ – periodic with respect to the variable $t$ and $C^2$ smooth with respect to $z$ and $\mu$.

2.2. Impact conditions

Suppose that the system (1) is defined for $z \in \Lambda$ and if $x_1 = 0$ the following impact conditions takes place. Let $\mu_- \leq 0 \leq \mu_+$ be a segment, $r : [\mu_-, \mu_+] \times \partial\Lambda \to [0, 1]$ be a $C^2$ smooth function.

**Condition 1.**

1. If the solution
$$z(t) = \mathrm{col}(z_1(t), \ldots, z_n(t)) =$$
$$\mathrm{col}(x_1(t), y_1(t), \ldots, x_n(t), y_n(t))$$

   of the system (1) is such that
$$x_1(t_0 - 0) = 0, \qquad y_1(t_0 - 0) \leq 0,$$

   then
$$x(t_0 + 0) = x(t_0 - 0),$$
$$y_1(t_0 + 0) = -r(\mu) y_1(t_0 - 0),$$
$$y_{tan}(t_0 + 0) = y_{tan}(t_0 - 0)$$

   where $y_\alpha(t_0 - 0)$ and $y_\alpha(t_0 + 0)$ are velocities before and after impact respectively, $\alpha$ is $1$ or $tan$.

2. Let the solution $z(t)$ of the system (1) be such that $x_1(t_0) = 0$, $y_1(t_0) = 0$ for a certain instant $t_0$. Let $\zeta(t)$ be the solution of the system
$$\dot{x}_k = y_k; \qquad \dot{y}_k = f_k(t, 0, 0, x_{tan}, y_{tan}, \mu),$$
$$k = 2, \ldots, n.$$

   with initial conditions $\zeta(t_0) = z_{tan}(t_0)$, and a segment
$$I = [t_0, t_1]$$

   be such that
$$f_1(t, 0, 0, \zeta(t), \mu) \leq 0$$

for all $t_1 \in I$. Then
$$x_1(t) = 0, \quad y_1(t) = 0$$

and $z_{tan}(t) = \zeta(t)$ for any $t \in I$.

**Remark 1**. So, we suppose, that $x_1$ is an impacting variable, the delimiter is homogeneous (the restitution coefficient does not depend neither on the point $x_{tan}$ nor on the velocity $y$ of the impact) and slippery (the tangent component of the velocity does not change during the impact).

The item 1 of the Condition 1 describes the instantaneous reflection, corresponding to a non-degenerate impact (with a positive value of the normal velocity), and the item 2 describes the sliding along the delimiter, corresponding to a zero-velocity impact. This part of the Condition 1 forbids the zero-velocity penetration through the delimiter which is a nonsense from the physical point of view.

We study the following model of vibro-impact system:

$$\dot{z} = F(t, z, \mu)$$
$$\text{while } x_1(t) > 0; \qquad (2)$$
$$\text{the Condition 1 is applied if } x_1(t - 0) = 0.$$

**Remark 2.** One can keep in mind the following simple case of the system (2):

$$\dot{x} = y; \qquad \dot{y} = -x - 0.1y + \sin t, \qquad x, y \in \mathbb{R}$$

while $x > 0$ and $y(t + 0) = -y(t - 0)$ if $x(t) = 0$. The mechanical system, corresponding to this model, is shown at the Figure 2. Here 1 is a mass, 2 is a spring, 3 is a delimiter and 4 is a damping element.

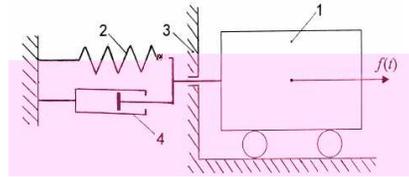

Fig. 2.

The following statement is evident due to the implicit function theorem.

**Lemma 1.** Let $z(t) = \mathrm{col}(x_1(t), y_1(t), \ldots, x_n(t), y_n(t))$ be the solution of the system (2) for $\mu = \mu_0$ and the initial conditions $z(t_0) = z^0 = \mathrm{col}(x_1^0, y_1^0, \ldots, x_n^0, y_n^0)$ where $x_1^0 \neq 0$. Suppose that this solution is defined on the segment $[t_-, t_+] \ni t_0$. Assume that there are exactly $N$ zeros

$$t_- < \tau_1^0 < \ldots < \tau_N^0 < t_+$$

of the function $x_1(t)$ over the segment $[t_-, t_+]$ and

$$y_1(\tau_j^0 - 0) \neq 0, \quad (j = 1, \ldots, N).$$

Then for any $\varepsilon > 0$ there exists a neighborhood $U$ of the point $\mathrm{col}(z^0, \mu_0)$ such that for any fixed

$$t \in \mathcal{J} = [t_-, t_+] \setminus \bigcup_{k=1}^{N}(\tau_k^0 - \varepsilon, \tau_k^0 + \varepsilon)$$

the mapping $z(t, t_1, z^1, \mu)$ is $C^2$ smooth. For these solutions there are exactly $N$ impact instants

$$\tau_j(t_1, z^1, \mu_1), \quad j = 1, \ldots, N$$

over the segment $[t_-, t_+]$. These instants and corresponding velocities

$$Y_j = -y_1(\tau_j(t_1, z^1, \mu_1) - 0, t_1, z^1, \mu_1)$$

$C^2$ smoothly depend on $t_1, z^1$ and $\mu_1$.

So, the instants and velocities of non-degenerate impacts smoothly depend on the initial data and parameters.

## 3. Grazing in the Newtonian model

Let $\mu^* \in (0, \mu_+)$. Suppose the following condition is satisfied.

**Condition 2.** (Fig. 3) *There exists a continuous family of $T$ – periodic solutions*

$$\varphi(t, \mu) = \mathrm{col}(\varphi_1^x(t, \mu), \varphi_1^y(t, \mu), \ldots, \varphi_n^x(t, \mu), \varphi_n^y(t, \mu)),$$

$t \in \mathbb{R}$, $\mu \in [0, \mu^*)$ *of the system (2) with the following properties.*

1. *For any $\mu \geq 0$ there exist exactly $N+1$ distinct zeros $\tau_0(\mu), \ldots, \tau_N(\mu)$ of the component $\varphi_1^x(t, \mu)$ over the period $[0, T)$.*

2. *The velocities $y_{0k}(\mu) = -\varphi_1^y(\tau_k(\mu) - 0, \mu)$ are such that*

$$\begin{aligned} y_{00}(\mu) &> 0 \quad \text{for all} \quad \mu > 0, \quad y_{00}(0) = 0, \\ f_1(\tau_0(0), 0, 0, \varphi_{tan}(\tau_0(0), 0), 0) &= \phi_0 > 0, \\ y_{0k}(\mu) &> 0, \quad \forall \mu \in [0, \mu^*), \quad k = 1, \ldots, N. \end{aligned} \quad (3)$$

*Here*

$$\varphi_{tan}(t, \mu) = (\varphi_2^x(t, \mu), \varphi_2^y(t, \mu), \ldots, \varphi_n^x(t, \mu), \varphi_n^y(t, \mu)).$$

3. *The instants $\tau_k(\mu)$ and the velocities $y_{0k}(\mu)$ continuously depend on $\mu \in [0, \mu^*)$.*

**Remark 3.** So we consider the family of periodic motions, corresponding to $N$ non-degenerate impacts over the period and exactly one impact, such that the normal component of the corresponding velocity vanishes for $\mu = 0$. This family may persist for negative values of $\mu$ (the so called *continuous grazing*, Figure 3) or disappear (a socalled *discontinuous grazing*). In the rest of this section we write down the conditions, sufficient for both of these bifurcations. A simple example of the system with a grazing family of periodic solutions (Eq. (27)) is studied at the Section 5.

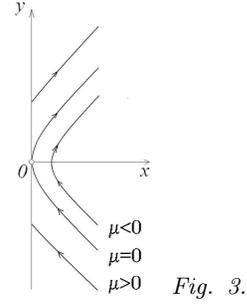

Fig. 3.

We suppose without loss of generality that $\tau_0(\mu) \equiv 0$, $\varphi_{tan}(0, \mu) \equiv 0$. This may be obtained by the transformation $s = t - \tau_0(\mu)$, $\zeta = z - \mathrm{col}(0, 0, \varphi_{tan}(0, \mu))$. Define

$$\theta_0 = \frac{1}{2} \min \left( \min_{\mu \in [0, \mu^*)} \tau_1(\mu), T - \max_{\mu \in [0, \mu^*)} \tau_N(\mu) \right).$$

Fix small positive $\theta < \theta_0$ and $\mu \in (0, \mu^*)$ consider the shift mapping for the system (2), given by the formula

$$S_{\mu, \theta}(z^0) = z(T - \theta + 0, -\theta, z^0, \mu).$$

For small positive $\mu$ and $\theta$ the mapping $S_{\mu, \theta}$ is $C^2$ smooth in a neighborhood of the point $z_{\mu, \theta} = \varphi(-\theta, \mu)$. Then

$$\lim_{\mu, \theta \to 0+} z_{\mu, \theta} = 0.$$

The existence of a grazing family of periodic solutions implies that the boundary problem for the system (2) with boundary conditions

$$x_1(\tau_0) = 0, \quad z(\tau_0 + 0) = P(\mu)z(\tau_0 + T - 0) \quad (4)$$

is solvable for small $\mu \geq 0$. Here $\tau_0$ is an unknown parameter,
$$P(\mu) = \begin{pmatrix} 1 & 0 & 0 \\ 0 & -r(\mu) & 0 \\ 0 & 0 & E_{2n-2} \end{pmatrix},$$
where $E_{2n-2}$ is the unit matrix of the corresponding size. Suppose that the considered boundary problem is solvable for $\mu = 0$. Without loss of generality, we assume that the corresponding values of $\tau_0$ and $z(\tau_0)$ are zeros.

Denote
$$\begin{aligned} A = (a_{ij}) &= (A_1, \ldots, A_{2n}) = (A_1, A_2, A_{tan}) = \\ &\lim_{\theta \to 0+} \lim_{\mu \to 0+} \frac{\partial z}{\partial z_0}(T - \theta - 0, \theta + 0, z_0, \mu)|_{z_0 = z_{\mu,\theta}}; \\ B = \mathrm{col}(b_1, \ldots, b_{2n}) &= \\ &\lim_{\theta \to 0+} \lim_{\mu \to 0+} \frac{\partial z}{\partial \mu}(T + \theta - 0, \theta + 0, z_0, \mu)|_{z_0 = z_{\mu,\theta}}; \\ F_0 &= \mathrm{col}(0, f_{01}, 0, f_2(0,0,0), \ldots, 0, f_n(0,0,0)). \end{aligned} \quad (5)$$

Now we find the conditions, sufficient for existence of a family $(\tau_0, z_0) = (\mathrm{col}(\tau_0(\mu), 0, y_{00}(\mu), z_{0,tan}(\mu)))$ solving the boundary problem for the system (2) with the boundary conditions (4) for small positive $\mu$. It follows from the definition of grazing that the inequality $y_{00}(\mu) > 0$ must be true for positive values of $\mu$.

Note that $b_1 = 0$ in the considered case. This statement follows from the Taylor formula for $x_1$.

Let $e_k$ be the $k$-th coordinate vector, $k = 1, \ldots, 2n$, $E_{tan} = \mathrm{col}(e_3, \ldots, e_{2n})$. Due to the implicit function theorem the existence of smooth functions $\tau_0(\mu)$, $y_{00}(\mu)$ and $z_{0,tan}(\mu)$ in the neighborhood of $\mu = 0$ can be provided by the condition
$$\det L_y \neq 0;$$
$$L_y = \left((E - A)F_0, A_2 + \frac{1}{r(0)}e_2, A_{tan} - E_{tan}\right). \quad (6)$$
Here $F_0 = F(0,0,0)$. The derivatives of $\tau_0(\mu)$, $y_{00}(\mu)$ and $z_{0,tan}(\mu)$ satisfy the system
$$L_y \mathrm{col}(\tau'(0), y'_{00}(0), z'_{0,tan}(0)) = -B.$$
So, to provide the inequality $y_{00}(\mu) > 0$ for positive values of $\mu$, we must have
$$[L_y^{-1} B]_2 < 0. \quad (7)$$
For the single-degree-of-freedom case the conditions (6) and (7) may be simplified. The matrix $L_y$ takes the form
$$\begin{pmatrix} -a_{12} f_{01} & a_{12} \\ (1 - a_{22}) f_{01} & a_{22} + 1/r(0) \end{pmatrix}.$$

The determinant of this matrix equals to
$$\Delta_y = -a_{12} f_{01}(1 + 1/r(0)).$$
So the condition (6) may be rewritten in the form $a_{12} \neq 0$. The inverse matrix $L_y^{-1}$ takes the form
$$\frac{1}{\Delta_y}\begin{pmatrix} a_{22} + 1/r(0) & -a_{12} \\ (a_{22} - 1) f_{01} & -a_{12} f_{01} \end{pmatrix}$$
and the condition (7) takes the form
$$b_2 < 0. \quad (8)$$

Let us discuss, when the family of periodic solutions, corresponding to grazing, persists for after-grazing (negative) values of $\mu$. In this case, the boundary problem for the system (2) with the conditions
$$y_1(\tau_0) = 0, \qquad z(\tau_0) = z(\tau_0 + T), \qquad x_1(\tau_0) \geq 0 \quad (9)$$
must be solvable.

Without loss of generality, we may suppose that it is the initial conditions $\tau_0 = 0$, $z(0) = 0$, which correspond to the solution of the boundary problem (9) for the system (2), $\mu = 0$. Then the existence of the smooth functions $\tau(\mu)$, $x_1(\mu)$, $z_{0,tan}(\mu)$, giving a solution of the considered boundary problem may be provided by the condition
$$\det L_x \neq 0; \qquad L_x = (A_1 - e_1, (E - A)F_0, A_{tan} - E_{tan}) \quad (10)$$
Since $f_{01} > 0$, the condition (10) may be rewritten in the simpler form $\det(A - E) \neq 0$.

To provide the condition $x_{01}(\mu) \geq 0$ for negative values of $\mu$ it suffices to suppose that
$$[L_x^{-1} B]_1 > 0, \quad (11)$$
where $[L_x^{-1} B]_1$ is the first element of the vector $L_x^{-1} B$.

In the s.d.f. case the condition (11) may be rewritten in a simpler form. If $n = 1$,
$$L_x = \begin{pmatrix} a_{11} - 1 & -a_{12} f_{01} \\ a_{21} & (1 - a_{22}) f_{01} \end{pmatrix}.$$
The determinant of this matrix equals to
$$\Delta_x = -\det(A - E) f_{01}$$
and
$$L_x^{-1} = \frac{1}{\Delta_x}\begin{pmatrix} (1 - a_{22}) f_{01} & a_{12} f_{01} \\ -a_{21} & a_{11} - 1 \end{pmatrix}$$

and the condition (11) is equivalent to the following one:

$$a_{12} b_2 \det(A - E) < 0. \qquad (12)$$

If the conditions (8) and (12) are satisfied simultaneously, $a_{12} \det(A - E) > 0$.

Finally, the following two scenarios are possible (the similar results for the square root type mapping are given in [11]).

1. A smooth family $\varphi(t, \mu)$ of periodic solutions exists for $\mu \in (\mu_-, \mu_+) \ni 0$. The number of impacts increases by 1 as the parameter passes through the bifurcation value. This case is called *continuous grazing bifurcation*. This bifurcation takes place if (7) and (11) are satisfied. For the s.d.f. case this corresponds to the case $b_2 < 0$, $a_{12} \det(A - E) > 0$.

2. A smooth family $\varphi(t, \mu)$ of periodic solutions exists for $\mu \in [0, \mu_+)$. These solutions disappear as $\mu$ passes through the bifurcation value. This case is called *discontinuous grazing bifurcation* and corresponds to inequalities (7) and

$$[L_x^{-1} B]_1 < 0. \qquad (13)$$

For the s.d.f. case the inequality (13) corresponds to the case $b_2 < 0$, $a_{12} \det(A - E) < 0$.

**Remark 4.** We mention the following classification of periodic solutions in vibro-impact systems, introduced by Peterka [1]. We say that the periodic solution $x(t)$ of a $T-$periodic vibro-impact system is of the $(m, n)$ type (sometimes the notation $m/n$ is used) if it is of the period $nT$ with exactly $m$ impacts over the period. Some bifurcations imply the changes in the Peterka pattern. For example, the non-degenarate period doubling corresponds to the transition $(m, n) \to (2m, 2n)$. In this sense the continuous grazing bifurcation corresponds to the transition $(m, n) \to (m \pm 1, n)$.

The bifurcation theory approach to this problem has been presented by A. P. Ivanov [48] for the s.d.f. case. The following classification has been given:

1. If $a_{12} < 0$ ($z_{12}$ in the notation of the quoted article) then for a positive value of the parameter $\mu$ one of the multipliers of the given *stable* periodic orbit turns to unity and this orbit collides with the unstable one and disappears, i.e. the saddle-node bifurcation scenario takes place.

2. If $a_{12} > 0$ then for a positive value of the parameter $\mu$ one of the multipliers of the given stable periodic orbit turns to $-1$ and the period doubling bifurcation takes place. In this case several periodic solutions of different periods may coexist.

Note that the validity of (8) depends on the choice of sign of $\mu$ and, if the considered periodic orbit is stable, we have $\det(A - E) > 0$. Hence, the existence of a continuous grazing family depends on the coefficient $a_{12}$ only.

## 4. Nonhyperbolic chaos in a neighborhood of the grazing bifurcation

Sometimes it may be difficult to apply the Devaney's definition of the chaotic invariant set to a multi-dimensional VIS. The main trouble is to check the hyperbolicity of the obtained invariant set. As we are going to show, one of Lyapunov exponents of the periodic solution in the neighborhood of grazing is large and positive, another one is large and negative, and nothing can be said about other ones, except they are relatively small with respect to the first two. Moreover, these exponents are very sensible to the parameters of the system and the choice of the mathematical model.

On the other hand, there are naturally defined stable and unstable manifolds of near-grazing periodic solutions (both of the dimension 1), corresponding to and the central manifold (see the Subsection 4.1 for definition) of the codimension 2. From the mechanical point of view, this means that the normal component $x_1$ of the vector $x$ is "fast" and the tangent component $x_{tan}$ is "slow".

### 4.1. Invariant manifolds

Let $U$ be a domain in the Euclidean space $\mathbb{R}^m$. Consider a diffeomorphism $S \in C^1(U \to V \subset \mathbb{R}^m)$. Suppose there exists an fixed point $x^* \in U$ of the mapping $S$. We can assume without loss of generality that $x^* = 0$.

Consider the constants

$$0 < \lambda_1 \leq \mu_1 < \lambda_2 \leq \mu_2 < \lambda_3 \leq \mu_3,$$

$\mu_1 < 1$, $\lambda_3 > 1$ and a decomposition $E^s \oplus E^u \oplus E^c = \mathbb{R}^m$, $DS(0)(E^\sigma) = E^\sigma$ for all $\sigma \in \{s, u, c\}$, such that the following conditions are satisfied. Here $\mu_1$ is the degree of contraction over the stable space and $\lambda_3$ is the degree of expansion over the unstable space.

1. Those eigenvalues $\lambda$ of the matrix $DS(0)$, which correspond to the vectors of the space $E^s$, satisfy the estimate $\lambda_1 \leq |\lambda| \leq \mu_1$.

2. Those eigenvalues $\lambda$, which correspond to the vectors of the space $E^c$, satisfy the estimate $\lambda_2 \leq |\lambda| \leq \mu_2$.

3. Those eigenvalues $\lambda$, which correspond to the vectors of the space $E^u$, satisfy the estimate $\lambda_3 \leq |\lambda| \leq \mu_3$.

The spaces $E^s$, $E^u$ and $E^c$ are called stable, unstable and central space respectively. Denote their dimensions by $m^s$, $m^u$ and $m^c$. Also, we consider the spaces

$$E^{cs} = E^s \oplus E^c, \quad E^{cu} = E^u \oplus E^c.$$

Following the notations of [44], we call them central stable and central unstable spaces. Let

$$m^{cs} = m^s + m^c, \quad m^{cu} = m^u + m^c.$$

The following result is a corollary of the reduction principle [46] (see also [43] and [47]) and an analogue of the well-known Perron theorem for a hyperbolic fixed point.

**Theorem 1.** *Let $S \in C^1(U \to \mathbb{R}^m)$, $0 \in U$ be a fixed point of the mapping $S$. Then there exist $C^1$ - smooth embeddings $b^s : \mathbb{R}^{m^s} \to \mathbb{R}^m$, $b^u : \mathbb{R}^{m^u} \to \mathbb{R}^m$, $b^{cs} : \mathbb{R}^{m^{cs}} \to \mathbb{R}^m$, $b^{cu} : \mathbb{R}^{m^{cu}} \to \mathbb{R}^m$, such that the following statements hold true.*

1. $b^s(0) = b^u(0) = b^{cs}(0) = b^{cu}(0) = 0$.

2. *Let $b^\sigma(\mathbb{R}^{m^\sigma}) = W^\sigma_{loc}$, $\sigma \in \{s, u, cs, cu\}$. These manifolds (we call them local stable, unstable, central stable and central unstable respectively) are locally invariant in a neighborhood of $0$.*

3. *The tangent spaces $T_0(W^\sigma_{loc})$ to $W^\sigma_{loc}(0)$ at the origin coincide with the spaces $E^\sigma$, $\sigma \in \{s, u, cs, cu\}$.*

4. *For any $\varepsilon > 0$ there exists positive numbers $a$, $\rho$ and $\delta$, such that*

$$|S^k(z)| \leq (a+\varepsilon)(\mu_1 + \varepsilon)^k |z|$$
*if* $z \in W^s_{loc}$, $|z| < \delta$, $k \in \mathbb{N}$;
$$|S^k(z)| \leq (a+\varepsilon)(\lambda_3 - \varepsilon)^{-k} |z|$$
*if* $z \in W^u_{loc}$, $|z| < \delta$, $-k \in \mathbb{N}$;
$$|S^k(z)| \leq (a+\varepsilon)(\mu_2 + \varepsilon)^k |z|$$
*if* $z \in W^{cs}_{loc}$, $|z| < \delta$, $k \in \mathbb{N}$
*while* $|S^k(z)| < \rho$;
$$|S^k(z)| \leq (a+\varepsilon)(\lambda_2 - \varepsilon)^{-k} |z|$$
*if* $z \in W^{cu}_{loc}$, $|z| < \delta$, $-k \in \mathbb{N}$
*while* $|S^k(z)| < \rho$.

5. *There exist numbers $\alpha > 0$ and $\varrho > 0$, such that for any $z \in B_\alpha(0) \setminus W^s_{loc}(0)$ there is a natural number $k_+$ such that $|S^{k_+}(z)|(\mu_1 - \varepsilon)^{k_+} \geq \varrho$, and for any $x \in B_\delta(0) \setminus W^u_{loc}(p)$ there is an integer $k_- < 0$ such that $|S^{k_-}(x)|(\lambda_3 - \varepsilon)^{k_-} \geq \varrho$.*

**Remark 5.** The stable and the unstable manifolds correspond to "fast" coordinates and the central manifold corresponds to "slow" coordinate.

We extend the obtained local manifolds to global ones, supposing

$$W^\sigma(0) = \{S^k(z) : k \in \mathbb{Z}, z \in W^\sigma_{loc}\},$$

$\sigma \in \{s, u, cs, cu\}$.

The behavior of iterations $S^k(p)$ of a point $p$ is the following:

1. $S^k(p)$ tends to zero exponentially as $k \to +\infty$ if

$$p \in W^s(0);$$

2. $S^k(p)$ tends to zero exponentially as $k \to -\infty$ if

$$p \in W^u(0);$$

3. $S^k(p)$ does not tend to zero too quickly (this is controlled by the parameter $\mu_2$) as $k \to +\infty$ if

$$p \in W^{cu}(0);$$

4. $S^k(p)$ does not tend to zero too quickly (this is controlled by the parameter $\lambda_2$) as $k \to -\infty$ if

$$p \in W^{cs}(0).$$

*4.2. The separatrix*

Let us consider the system (2), assuming the Condition 2 is satisfied, i.e. there is a grazing family of periodic solutions. Fix $\mu, \theta > 0$. Denote

$$\Gamma_{\mu,\theta} = \{z^0 \in \Lambda : \exists t_1 \in [-T, T] : z_1(t_1, -\theta, z^0, \mu) = 0\}$$

(Fig. 4).

**Lemma 2.** *There exists a neighborhood $U_0$ of zero such that if the positive parameters $\mu$ and $\theta$ are small enough, the set $\Gamma_{\mu,\theta} \bigcap U_0$ is a surface of the dimension $2n - 1$, which is*

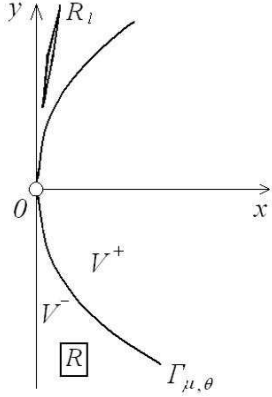

Fig. 4.

the graph of the $C^2$ smooth function $x_1 = \gamma_{\mu,\theta}(x_{tan}, y)$. Moreover,

$$\gamma_{\mu,\theta}(z) = y_1^2 \left( \frac{1}{f_1(-\theta, z_{\mu,\theta}, \mu)} + \widetilde{\gamma_{\mu,\theta}}(z) \right), \quad (14)$$

where $\widetilde{\gamma_{\mu,\theta}}$ is a $C^2$ smooth function such that $\widetilde{\gamma_{\mu,\theta}}(0) = 0$.

**Remark 6.** This surface $\Gamma_{\mu,\theta}$ plays an important role in our later reasonings. We show that corresponds to the points of non-smoothness of the Poincaré mapping and, at the same time, to the solutions with a degenerate impact. The invariant manifolds bend, intersecting the surface $\Gamma_{\mu,\theta}$.

**Proof.** Take a point $\zeta \in \Gamma_{\mu,\theta}$. Let the instant $t_0$ be such that $z_1(t_0, -\theta, \zeta) = 0$, $s = t - t_0$,

$$z_1(t + 0, -\theta, \zeta) = \mathrm{col}(x_1(t), y_1(t)).$$

Let us show that if $t_0$ is close enough to $-\theta$, we may take $s_0 \geq |t_0 + \theta|$ so that the function $x_1(t_0 + s)$ does not have zeros on $[-s_0, s_0]$, except $s = 0$. Otherwise, there exists a sequence $t_0^k \to -\theta$ (suppose without loss of generality, that $t_0^k > -\theta$ and the sequence decreases), a sequence $t_1^k \in [-\theta, t_0^k)$ and one, consisting of solutions, uniformly bounded on the segment $[-\theta, t_0^1]$:

$$z^k(t) = \mathrm{col}(z_1^k(t), \ldots, z_n^k(t)) = \\ \mathrm{col}(x_1^k(t), y_1^k(t), \ldots, x_n^k(t), y_n^k(t))$$

such that $z_1^k(t_0^k) = 0$, $x_1^k(t_1^k) = 0$ (Fig. 5). Also, there exist instants $t_2^k \in (t_1^k, t_0^k)$, such that $\dot{x}_1^k(t_2^k) = 0$ and instants $t_3^k \in (t_2^k, t_0^k)$ such that $\ddot{x}_1^k(t_3^k) = 0$. Moreover, $t_3^k \to -\theta$, $x_1^k(t_3^k) \to 0$, $\dot{x}_1^k(t_3^k) \to 0$. Then

$$\ddot{x}_1(t_3^k) \to f_1(-\theta, 0, \mu) = 0.$$

This contradicts to (3).

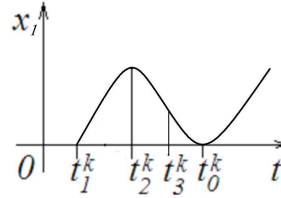

Fig. 5.

Then for all $s \in [-s_0, s_0]$ the function $x_1(t_0 + s)$ can be represented as series

$$x_1(t_0 + s) = X_2 s^2 + X_3 s^3 + \ldots \quad (15)$$

Differentiating (15), we obtain that

$$\dot{x}_1(t_0 + s) = 2X_2 s + 3X_3 s^2 + \ldots.$$

On the other hand,

$$X_2 = \ddot{x}_1(t_0 + 0)/2 \to f_1(-\theta, z_{\mu,\theta}, \mu)$$

as $t_0 \to -\theta$. Then

$$x_1(-\theta) = f_1(-\theta, z_{\mu,\theta}, \mu)(t_0 + \theta)^2(1 + o(1))/2;$$

$y_1(-\theta) = f_1(-\theta, z_{\mu,\theta}, \mu)(t_0 + \theta)(1 + o(1))$. Since $y_1 = \dot{x}_1$, the formula (14) is true. ∎

Take a small parameter $\varsigma > 0$ such that the sets

$$V = \{z \in \Lambda : \|z - z_{\mu,\theta}\| \leq \varsigma\} \subset U_0, \\ \{(x, y) \in V : x_1 < \gamma_{\mu,\theta}(x_{tan}, y)\}, \\ \{(x, y) \in V : x_1 > \gamma_{\mu,\theta}(x_{tan}, y)\}$$

are correctly defined and nonempty.

Consider the matrix $A$, defined by the formula (5). Let $\Delta_0 = \det A$. Denote the elements of the matrix $A$ by $a_{ij}$ and ones of the matrix $A^{-1}$ by $\alpha_{ij}$. Denote the columns of matrices $A$ and $A^2$ by $A_j$ and $A_j^2$ respectively, the strings of the matrix $A^{-1}$ by $\mathcal{A}_j$. Later on, we assume that

$$a_{12} > 0, \quad \sum_{j=1}^{2n} a_{1k} a_{k2} < -a_{12}. \quad (16)$$

As we show later, these conditions allow us to control the bent of invariant manifolds, intersecting the surface $\Gamma_{\mu,\theta}$.

### 4.3. The Jacobi matrix

Note that all the mappings $S_{\mu,\theta}$, corresponding to the same value of $\mu$, are conjugated. Fix a number $\mu > 0$

and a solution $z^0(t) = \text{col}(x_1^0(t), y_1^0(t), \ldots, x_n^0(t), y_n^0(t))$ of the corresponding system with an impact at $t_0$. Suppose that the corresponding normal velocity $Y_{01} = -y_1^0(t_0 - 0)$ is nonzero. We consider $Y_{01}$ as a small parameter. Denote $Z_{0,tan} = z_{tan}^0(t_0 - 0)$. Fix a positive number $s_0$ and consider the mapping

$$\Sigma(\zeta) = z(t_0 + s_0, t_0 - s_0, \zeta, \mu),$$

defined in a neighborhood of the point $\zeta^0 = z^0(t_0 - s_0)$. Here we assume that the point $\zeta^0$ and the parameter $s_0$ are chosen so that there exists a neighborhood $\Omega \ni \zeta^0$ such that any solution $z(t) = z(t, t_0 - s_0, z_-, \mu)$ ($z_- \in \Omega$) impacts once over the segment $[t_0 - s_0, t_0 + s_0]$. Denote the corresponding instant by $t_1 = t_1(z_-)$ and the normal velocity of the impact by $Y_1 = Y_1(z_-)$. Let $Z_{tan} = z_{tan}(t_1)$ be the tangent component of the solution $z(t)$ at the impact instant. Take the numbers $s_\pm = s_\pm(z_-)$ so that

$$t_0 \pm s_0 = t_1(z_-) \pm s_\pm(z_-)$$

for all $z_- \in \Omega$. The mapping $\Sigma$ is smooth in the neighborhood of the point $\zeta_0$, let us estimate the Jacobi matrix $D\Sigma(\zeta_0)$. Denote

$$z_+ = z(t_0 + s_0) = z(t_0 + s_0, t_1 + 0, 0, rY_1, Z_{tan}, \mu),$$
$$x_+ = x(t_0 + s_0) = x(t_0 + s_0, t_1 + 0, 0, rY_1, Z_{tan}, \mu),$$
$$y_+ = z(t_0 + s_0) = y(t_0 + s_0, t_1 + 0, 0, rY_1, Z_{tan}, \mu),$$
$$x_- = x(t_0 - s_0) = x(t_0 - s_0, t_1 - 0, 0, -Y_1, Z_{tan}, \mu),$$
$$y_- = y(t_0 - s_0) = y(t_0 - s_0, t_1 - 0, 0, -Y_1, Z_{tan}, \mu).$$

Similarly, we define the numbers

$$x_{1,\pm}, \quad y_{1,\pm}, \quad x_{\pm,tan}, \quad y_{\pm,tan}.$$

Consider the Taylor formula for $z_\pm$:

$$x_{1,-} = Y_1 s_- + f_1(t_1, 0, -Y_1, Z_{tan}, \mu) s_-^2/2 +$$
$$\rho_1(s_-, t_1, Y_1, Z_{tan}, \mu) s_-^3;$$
$$y_{1,-} = -Y_1 - f_1(t_1, 0, -Y_1, Z_{tan}, \mu) s_- +$$
$$\rho_2(s_-, t_1, Y_1, Z_{tan}, \mu) s_-^2;$$
$$x_{-,tan} = x_{tan}(t_1) - y_{tan}(t_1 - 0) s_- +$$
$$f_{tan}(t_1, 0, -Y_1, Z_{tan}, \mu) s_-^2/2 +$$
$$\rho_3(s_-, t_1, Y_1, Z_{tan}, \mu) s_-^3;$$
$$y_{tan-} = y_{tan}(t_1 - 0) - f_{tan}(t_1, 0, -Y_1, Z_{tan}, \mu) s_- +$$
$$\rho_4(s_-, t_1, Y_1, Z_{tan}, \mu) s_-^2;$$
$$x_{1,+} = rY_1 s_+ + f_1(t_1, 0, rY_1, Z_{tan}, \mu) s_+^2/2 +$$
$$\rho_5(s_+, t_1, Y_1, Z_{tan}, \mu) s_+^3;$$
$$y_{1,+} = rY_1 + f_1(t_1, 0, rY_1, Z_{tan}, \mu) s_+ +$$
$$\rho_6(s_+, t_1, Y_1, Z_{tan}, \mu) s_+^2;$$
$$x_{tan+} = x_{tan}(t_1) + y_{tan}(t_1 - 0) s_+ +$$
$$f_{tan}(t_1, 0, rY_1, Z_{tan}, \mu) s_+^2/2 +$$
$$\rho_7(s_+, t_1, Y_1, Z_{tan}, \mu) s_+^3;$$
$$y_{tan+} = y_{tan}(t_1 - 0) + f_{tan}(t_1, 0, rY_1, Z_{tan}, \mu) s_+ +$$
$$\rho_8(s_+, t_1, Y_1, Z_{tan}, \mu) s_+^2.$$
(17)

Here all functions, denoted by the letter $\rho$ with different indices, are $C^2$ smooth with respect to all arguments except $s_\pm$. Denote

$$f_{0k+} = f_k(t_0, 0, rY_{01}, Z_{0,tan}, \mu),$$
$$f_{0k-} = f_k(t_0, 0, -Y_{01}, Z_{0,tan}, \mu).$$

It follows from (17) that

$$\left.\frac{\partial z_+}{\partial(s_+, Y_1, Z_{tan})}\right|_{s_+=0,\, Y_1=Y_{01},\, Z_{tan}=Z_{0,tan}} =$$
$$\begin{pmatrix} rY_{01} & 0 & 0 \\ f_{01+} & r & 0 \\ X_+ & 0 & E_{2n-2} \end{pmatrix};$$
$$\left.\frac{\partial z_-}{\partial(s_-, Y_1, Z_{tan})}\right|_{s_-=0,\, Y_1=Y_{01},\, Z_{tan}=Z_{0,tan}} =$$
$$\begin{pmatrix} Y_{01} & 0 & 0 \\ -f_{01-} & -1 & 0 \\ X_- & 0 & E_{2n-2} \end{pmatrix}.$$

Here $E_{2n-2}$ is the unit matrix of the corresponding size, $X_+ = \text{col}(y_2(t_0 - 0), f_{02+}, \ldots, y_n(t_0 - 0), f_{0n+})$,

$$X_- = \text{col}(-y_2(t_0 - 0), -f_{02+}, \ldots, -y_n(t_0 - 0), -f_{0n+}).$$

Denote

$$f_k' = \left.\frac{\partial f_k(t_0, 0, y_1, Z_{0,tan}, \mu)}{\partial y_1}\right|_{y_1=0}.$$

Clearly, $ds_+/ds_- = -1$. Then, similarly to the results of the paper [49], we obtain

$$\mathcal{B} = \lim_{s_\pm \to 0} \frac{\partial z_+}{\partial z_-} = \begin{pmatrix} -r & 0 & 0 \\ \beta_{21} & -r & 0 \\ \mathcal{B}_{1,tan} & 0 & E_{2n-2} \end{pmatrix}.$$

Here $\mathcal{B}_{1,tan} = \mathrm{col}(\beta_{31}, \ldots \beta_{2n\,1})$,

$$\beta_{21} = -(f_{01+} + r f_{01-})/Y_{01} =$$
$$-(r+1)\phi_0(1 + O(Y_{01}))/Y_{01},$$
$$\beta_{2j-1\,1} = 0,$$
$$\beta_{2j\,1} = (f_{0j+} - f_{0j-})/Y_{01} = (r+1)f'_k + O(Y_{01}),$$
$$j = 1, \ldots, n.$$

Note that $\det \mathcal{B} = r^2$.

It is shown at the Figure 4 how a small rectangular domain $U$ of initial data is stretched while the time passes a fixed segment corresponding to low-velocity impacts of solutions ($R$ transfers to $R_1$).

*4.4. Lyapunov exponents*

In this subsection we study the Lyapunov exponents, corresponding to the Poincaré mapping at a near-grazing periodic point. We show that one of these exponent tend to $+\infty$ as $\mu \to 0$, another one tends to $-\infty$ and all others are relatively small, compared with the first two (the ratios tend to zero as $\mu \to 0$).

We estimate the larger and the smaller absolute value of eigenvalues of $D = DS_{\mu,\theta}(z_{\mu,\theta})$ and ones of small perturbations of this matrix. The mapping $S_{\mu,\theta}$ can be represented as the composition $S_{\mu,\theta} = S_{2,\mu,\theta} \circ S_{1,\mu,\theta}$, where

$$S_{2,\mu,\theta}(\zeta) = z(T - \theta + 0, \theta, \zeta, \mu).$$

The Jacobi matrix $A_{\mu,\theta} = DS_{2,\mu,\theta}(S_{1,\mu,\theta}(z_{\mu,\theta}))$ tends to $A$ as $\mu, \theta \to 0$. Note, that the elements of these matrices satisfy the inequalities similar to (16) provided $\mu$ and $\theta$ are sufficiently small.

The matrix $\mathcal{B}_{\mu,\theta} = DS_{1,\mu,\theta}(z_{\mu,\theta})$ is of the form

$$\lim_{s_\pm \to 0} \frac{\partial z_+}{\partial z_-} = \mathcal{B}$$

as $\mu$ and $\theta$ tend to 0.

Then

$$D = A_{\mu,\theta} \mathcal{B}_{\mu,\theta} =$$
$$(-(r+1)A_2\phi_0(1 + o(1))/Y_0, -rA_2(1 + o(1)),$$
$$A_3 + o(1), \ldots, A_{2n} + o(1)).$$

Consequently, if $a_{12} \neq 0$ and if $\mu$ is small, one of the eigenvalues of the matrix $D$ is $\lambda_+ = -(r+1)a_{12}\phi_0(1 + o(1))/Y_0$. The corresponding eigenvector $u_+$ equals to $A_2 + o(1)$. The eigenvalue $\lambda_+$ is of the multiplicity 1, the linear space, corresponding to other eigenvalues, tends to the hyperplane $\pi_1$, given by the condition $x_1 = 0$ as $\mu \to 0$. Since $a_{12} \neq 0$, the vector $u^+$ is out of $\pi_1$. Note that

$$\det D = \det A_{\mu,\theta} \det \mathcal{B}_{\mu,\theta} = (r^2 + o(1))\Delta_0.$$

The matrix $D^{-1}$ satisfies the following asymptotic estimate

$$D^{-1} = \frac{1}{r^2} \begin{pmatrix} -r\mathcal{A}_1 + o(1) \\ (r+1)\mathcal{A}_1(1 + o(1))/Y_0 \\ \mathcal{A}_3 + o(1) \\ \ldots \\ \mathcal{A}_{2n} + o(1) \end{pmatrix}.$$

It follows from the form of this matrix, one of the eigenvalues of the matrix $D^{-1}$ equals to

$$\lambda_-^{-1} = (r+1)\phi_0\alpha_{12}(1 + o(1))/(r^2 Y_0).$$

The corresponding eigenvector satisfies the asymptotical estimate $u_- = e_2 + o(1)$.

*4.5. Intersection of invariant manifolds*

Let

$$\lambda_+, \lambda_-, \lambda_3, \ldots, \lambda_{2n}$$

the eigenvalues of the matrix $DS_{\mu,\theta}(z_{\mu,\theta})$. Let

$$\Lambda_1 = \min(|\lambda_+|, |\lambda_-^{-1}|)/2, \qquad \Lambda_2 = 2 \max_{k=3,\ldots,2n}(|\lambda_k|, |\lambda_k^{-1}|).$$

It is only the first column of $D$, which is unbounded. Consequently, $\Lambda_1 > \Lambda_2$ provided $\mu$ is small enough. Let $E^s$ be the line, corresponding to the eigenvalue $\lambda_-$, $E^u$ be the one-dimensional space, corresponding to the eigenvalue $\lambda_+$, and $E^c$ be the space, corresponding to all other eigenvalues. The spaces $E^s$ and $E^s$ are linear hulls of the vector $E_1^s = e_2 + O(\theta)$ and the vector $E_1^u = A_2 + O(\theta)$ respectively.

Due to the Theorem 1, for the fixed point $z_{\mu,\theta}$ there exist the stable manifold $W^s$ the unstable manifold $W^u$ (both of the dimension 1), the central stable manifold $W^{cs}$ and the central unstable manifold $W^{cu}$ (both of the dimension $2n - 1$). The corresponding spaces $E^\sigma$ are tangent to these manifolds at the fixed point.

Now we show that, provided the assumptions of the Theorem 1 are satisfied there exists a non-hyperbolic analogue of a homoclinic point for the considered system.

**Lemma 3.** *There is $\mu_+ > 0$, so small that, provided*

$$\mu \in (0, \mu_+),$$

*the manifolds $W^{cs}$ and $W^u$ intersect transversally at a point $p \neq z_{\mu,\theta}$, i.e. there exists a neighborhood $V$ of the point $p$, such that the connected components of sets $W^{cs} \bigcap V$ and $W^u \bigcap V$, containing $p$, are smooth manifolds and intersect transversally at the point $p$.*

This intersection for the s.d.f. case is illustrated at the figure 6.

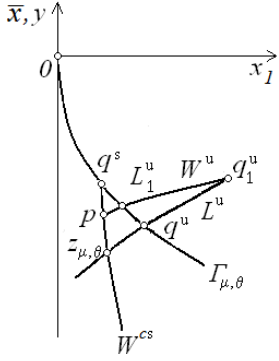

Fig. 6.

**Proof.** For small positive values of $\mu$ the manifold $W^u$ intersects transversally the surface $\Gamma_{\mu,\theta}$ at the point $q^u$. Let $q_1^u = S_{\mu,\theta}(q^u)$. Denote by $L^u$ the arc of the manifold $W^u$, connecting the points $q^u$ and $q_1^u$. Denote $L_1^u = S_{\mu,\theta}(L^u)$. The curve $W^u$ is not smooth at the point $q_1^u$. The tangent lines to $W^u$ at the points of the curve $L^u$ are close to the line directed along the vector $A_2$ i.e. the direction of the vector $A_2$ is the limit one to all tangent vectors of $W^u$ as $\mu, \theta \to 0$. The tangent lines to the arc $L_1^u$ are close to the linear hull of the vector $A_2^2$ for small $\mu$. Due to (16) the vectors $A_2$ and $A_2^2$ are directed to the different sides with respect to the hyperplane $\pi_1$ and of the space $E^{cs}$. Consequently, provided $\mu$ is small, the neighborhood $V$ can be chosen so that the arc $L_1^u$ and the manifold $W^{cs}$ intersect transversally. ∎.

*4.6. Straightening*

The main goal of this subsection is to introduce the appropriate coordinates in a neighborhood of a homoclinic point.

Due to the Theorem 1, there exists a diffeomorphism $h : U_{\mu,\theta} \to \Omega$ of the neighborhood $U_{\mu,\theta}$ of the point $z_{\mu,\theta}$ to a domain $\Omega$, containing the origin in $\mathbb{R}^{2n}$ and endowed with coordinates

$$\zeta = \mathrm{col}(\zeta^s, \zeta^u, \zeta^c) = \mathrm{col}(\zeta^s, \zeta^u, \zeta_1^c, \ldots, \zeta_{2n-2}^c)$$

and such that the following conditions are satisfied (here we suppose that $h$ is extended to a diffeomorphism of the whole Euclidean space).

1. $h(z_{\mu,\theta}) = 0$.

2. The Euclidean norms of all columns of the matrix $Dh(z_{\mu,\theta})$ equal to 1.

3. $h(W_{loc}^s) \subset \{\mathrm{col}(\zeta^s, 0, 0) : \zeta^s \in \mathbb{R}\} = O\zeta^s$,

$$h(W_{loc}^u) \subset \{0, 0, \zeta^u) : \zeta^u \in \mathbb{R}\} = O\zeta^u,$$

$$h(W_{loc}^{cs}) \subset \{\mathrm{col}(\zeta^s, \zeta^c, 0) : \zeta^s \in \mathbb{R}, \zeta^c \in \mathbb{R}^{2n-2}\},$$

$$h(W_{loc}^{cu}) \subset \{\mathrm{col}(0, \zeta^c, \zeta^u) : \zeta^u \in \mathbb{R}, \zeta^c \in \mathbb{R}^{2n-2}\}.$$

4. $Dh(z_{\mu,\theta})u_+ = e_1 = \mathrm{col}(1, 0, \ldots, 0)$.

5. $Dh(z_{\mu,\theta})u_- = e_2 = \mathrm{col}(0, 1, 0, \ldots, 0)$.

6. The point $q = h(p)$ corresponds to the coordinates $(\zeta_q^s, 0, \zeta_q^c)$, $\zeta_q^s > 0$. The segment, linking the point $q$ with the origin, is a subset of $h(W^{cs} \bigcap U)$.

7. Let $L_U^u$ be the connected component of $h(W^u \bigcap U)$, containing the point $q$. Then

$$h(L_U^u) \subset \{\mathrm{col}(\zeta_q^s, \zeta^u, \zeta_q^c) : \zeta^u \in \mathbb{R}\}.$$

Let $\hat{S}_{\mu,\theta} = h \circ S_{\mu,\theta} \circ h^{-1}$. Consider the point $q_1$ with the coordinates $(\zeta_q^s, 0, 0)$. Select the positive numbers $\varepsilon^s$, $\varepsilon^u$ and $\varepsilon^c$ and the neighborhoods

$$\begin{aligned} U_0 &= \{\zeta = \mathrm{col}(\zeta^s, \zeta^u, \zeta^c) : |\zeta^s| \leq \varepsilon^s, |\zeta^u| \leq \varepsilon^u, |\zeta^c| \leq \varepsilon^c\}; \\ U_1 &= \{\zeta = \mathrm{col}(\zeta^s, \zeta^u, \zeta^c) : |\zeta^s - \zeta_q^s| \leq \varepsilon^s, \\ &\quad |\zeta^u| \leq \varepsilon^u, |\zeta^c| \leq \varepsilon^c\} \end{aligned} \quad (18)$$

so that the following conditions are satisfied

1. $U_0, U_1 \subset U$.

2. $U_0 \bigcap U_1 = \emptyset$, i.e. $\varepsilon^s < \zeta_q^s/2$.

3. $U_1 \bigcap h(W^s)$ i.e. $\varepsilon^c > 2|\zeta_q^c|$.

4. Let $\Pi_1$ be the projector to the axis $O\zeta^s$. For any two points
$$\zeta_{1,2} = \operatorname{col}(\zeta_{1,2}^s, \zeta^u, \zeta^c) \in U$$
the following Lipschitz condition is true:
$$|\Pi_1 \hat{S}_{\mu,\theta}(\zeta_1) - \Pi_1 \hat{S}_{\mu,\theta}(\zeta_2)| \le |\zeta_1^s - \zeta_2^s|/2.$$

So, we have introduced two disjoint sets: a neighborhood $U_0$ of the fixed point and a neighborhood $U_1$ of the homoclinic point. The set $U_1$ is constructed so that it does intersect the local stable manifold. The reductions $\Pi_1 \hat{S}_{\mu,\theta}|_{U_i}$ must be uniformly contracting. This may be obtained since the distance between the homoclinic point and the fixed point tends to zero as the parameter $\mu$ tends to the bifurcation value (Fig. 7).

### 4.7. Admissible disks

Here we consider the set of possible approximations to central unstable manifolds and study their iterations.

We define an admissible disk in one of the domains $U_i$ as a subset of $U_i$, which is the graph of a $C^1$ – smooth function $\zeta^s = \eta(\zeta^u, \zeta^c)$, defined for $|\zeta^u| \le \varepsilon^u$, $|\zeta^c| \le \varepsilon^c$ and such that $\max|D\eta(\zeta^u, \zeta^c)| \le 1$. Let $\mathcal{Q}_i$ ($i=0,1$) be the sets of disks, admissible at $U_i$ endowed with the $C^1$ metric $\operatorname{dist}_1$. Denote $\mathcal{Q} = \mathcal{Q}_0 \bigcup \mathcal{Q}_1$.

Any admissible disk intersects transversally with the stable manifold. There exists an integer $m_1$ such that the image $\hat{S}_{\mu,\theta}^{m_1}(Q)$ of any admissible disk $Q$ does intersect the domain $U_0$. Moreover, the intersection contains an admissible curve $\varpi$ i.e. the graph of a $C^1$ – smooth function $(\zeta^s, \zeta^c) = \chi(\zeta^u)$, defined for $|\zeta^u| \le \varepsilon^u$ and such that
$$\max_{|\tau| \le \varepsilon^u} |\chi'(\zeta^u)| \le 1, \qquad \chi(0) = (\chi^s(0), 0)$$
i.e. any admissible curve intersects the local stable manifold. There exist positive numbers $\mu_0$ and $\varepsilon_0$ and the integer number $m_2$ such that if $\mu < \mu_0$, the image $\hat{S}_{\mu,\theta}^{m_2}(\varpi)$ of any admissible curve $\varpi$, satisfying the condition $|\chi(0)| \le \varepsilon_0$, contains an admissible curve $\varpi_1 \subset U_1$. Without loss of generality we may suppose $m_2$ to be so large that for any admissible disk $Q$ the intersection $\hat{S}_{\mu,\theta}^{m_1}(Q) \bigcap h(W^s)$ contains a point $\operatorname{col}(\zeta^s, 0, 0)$; $|\zeta^s| \le \varepsilon_0$.

**Remark 7.** Actually, we are going to present a structure, which approximates a foliation of central unstable manifolds. In a hyperbolic case we could take a disk, whose inclination to the stable manifold is not too small (this is what we call admissible disk). Then, due to the $\lambda$ – lemma the iterations of this disk provide local approximations to the unstable manifold. For non-hyperbolic cases the iterations of these admissible disks may shrink in some directions. This is why we develop a more complex technique, using the procedure of extension (see the next subsection).

Select $\varepsilon_1^c \in (0, \varepsilon^c]$ so that for any $m > m_1$ the images $\hat{S}_{\mu,\theta}^m(Q)$ of any admissible disk $Q$ contain disks, admissible at the more narrow domain
$$\widetilde{V}_0 = \{\zeta : |\zeta^s| \le \varepsilon^s, |\zeta^u| \le \varepsilon^u, |\zeta^c| \le \varepsilon_1^c\}.$$
Select the number $\varepsilon_2^c$ so that the image $\hat{S}_{\mu,\theta}^{m_2}(\widetilde{Q})$ of any disk $\widetilde{Q}$, admissible at the set $\widetilde{V}_0$, contains a disk, admissible at the domain
$$\{\zeta : |\zeta^s - \zeta_q^s| \le \varepsilon^s, |\zeta^u| \le \varepsilon^u, |\zeta^c - \zeta_q^c| \le \varepsilon_2^c\}$$
and one, admissible at the domain
$$\{\zeta : |\zeta^s| \le \varepsilon^s, |\zeta^u| \le \varepsilon^u, |\zeta^c| \le \varepsilon_2^c\}.$$

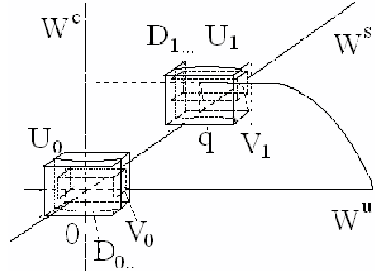

Fig. 7.

Denote $\delta^c = \min(\varepsilon_1^c, \varepsilon_2^c)$, $m_3 = m_1 + m_2$.

### 4.8. Extensions

We are going to show that the initial dynamical system engenders one at the set of admissible disks. Unfortunately, we cannot consider the iterations $S_{\mu,\theta}^k(D)$ of an admissible disk $D$ since these iterations may shrink. So, we need to "expand" these iterations on every step.

Consider the sets
$$V_0 = \{\zeta : |\zeta^s| \le \varepsilon^s, |\zeta^u| \le \varepsilon^u, |\zeta^c| \le \delta^c/2\},$$
$$V_1 = \{\zeta : |\zeta^s - \zeta_q^s| \le \varepsilon^s, |\zeta^u| \le \varepsilon^u, |\zeta^c - \zeta_q^s| \le \delta^c/2\},$$
$$\partial_0 = \{\zeta : |\zeta^s| \le \varepsilon^s, |\zeta^u| \le \varepsilon^u, |\zeta^c| = \delta^c/2\},$$
$$\partial_1 = \{\zeta : |\zeta^s - \zeta_q^s| \le \varepsilon^s, |\zeta^u| \le \varepsilon^u, |\zeta^c - \zeta_q^s| = \delta^c/2\}.$$

We can take $k$ so large that the following statements hold.

1. For all $i,j \in \{0,1\}$ the intersections

$$\hat{S}_{\mu,\theta}^k(U_i) \bigcap U_j$$

contain components $U_{ij}$ such that if $Q \in \mathcal{Q}_i$ then the intersection $\hat{S}_{\mu,\theta}^k(Q) \bigcap U_{ij}$ contains an admissible disk $\widetilde{G}_j(Q) \subset V_j \subset U_{ij}$.

2. For all $Q_{1,2} \in \mathcal{Q}$, $i = 0, 1$

$$\text{dist}_1(\widetilde{G}_i(Q_1), \widetilde{G}_i(Q_2)) \leq \frac{1}{2}\text{dist}_1(Q_1, Q_2). \quad (19)$$

3. The mappings $\widetilde{G}_i$ are one-to-one; moreover, the images of non-intersecting disks do not intersect any more.

Denote by $\mathcal{Q}^-$ the sets of the disks, admissible in corresponding domains $V_i$. Clearly, any disk $Q \in \mathcal{Q}$ contains a subdisk from the set $\mathcal{Q}^-$. Here we construct an embedding $J : \mathcal{Q}^- \to \mathcal{Q}$ such that

$$\text{dist}_1(J(Q_1), J(Q_2)) \leq \text{dist}_1(Q_1, Q_2) \quad (20)$$

for all $Q_{1,2} \in \mathcal{Q}$ and $J(Q) \bigcap V_i = Q \bigcap V_i$ for all $Q \in \mathcal{Q}^-$, $i = 0, 1$.

Fix a disk $\widetilde{Q} \in \mathcal{Q}^-$, given by the equation

$$\zeta^s = \eta_1(\zeta^u, \zeta^c).$$

Let $i \in \{0,1\}$ be such that $\widetilde{Q} \subset U_i$. The intersection of the disk $\widetilde{Q}$ and the surface $\partial_{V_i}$ is a manifold $\partial_Q$ of the dimension $2n-2$. Let us join every point of this manifold with a point of the set $\partial_{U_i}$ by a segment of the type

$$\{\zeta = \text{col}(\zeta^s, \zeta^u, \zeta^c) : \zeta^s = c_1, \zeta^u = c_2,$$
$$\zeta^c = \tau c_3 : \tau \in [\tau_-, \tau_+]\}$$

so that $(c_1, c_2, \tau_- c_3) \in \partial_Q$, $(c_1, c_2, \tau_+ c_3) \in \partial_{U_i}$ or vice versa. We obtain a disk in $U_i$, defined by formulae $\zeta^s = \eta_2(\zeta^u, \zeta^c)$. However, this disk is not smooth at the points of the set $\partial_Q$.

Let us smooth this disk in a standard way. Suppose, without loss of generality that $i = 0$. Consider a $C^\infty$ – smooth function $\phi : \mathbb{R} \to [0,1]$, satisfying the following conditions.

1. $\phi(s) \equiv 0$ if $s \leq 0$, $\phi(s) \equiv 1$ if $s \geq 1$.

2. $\phi'(s) \in (0,2)$ for all $s \in (0,1)$.

Let $B = \{z \in \mathbb{R}^{2n} : |z| \leq \varepsilon^c\}$. Define the function $\eta : [-\varepsilon^u, \varepsilon^u] \times B$ by formulae:

$$\eta(\zeta^u, \zeta^c) = \begin{bmatrix} \eta_1(\zeta^u, \zeta^c), & \text{if } |\zeta^c| \leq \delta^c/2; \\ (1 - \phi((2|\zeta^c| - \delta^c)/\delta^c))\eta_1(\zeta^u, \zeta^c) + \\ \phi((2|\zeta^c| - \delta^c)/\delta^c)\eta_2(\zeta^u, \zeta^c), \\ \text{if } |\zeta^c| \in [\delta^c/2, \delta^c]; \\ \eta_2(\zeta^u, \zeta^c), & \text{if } \quad |\zeta^c| \geq \delta^c. \end{bmatrix}$$

Clearly, the graph of this function $\eta$ is an admissible disk.

Now we fix an admissible disk $Q \subset U_0 \bigcup U_1$ and apply the discussed extension-smoothing procedure to the disks $\widetilde{G}_0(Q)$ and $\widetilde{G}_1(Q)$. Denote the obtained disks by $G_0(Q)$ and $G_1(Q)$. The mappings $G_0$ and $G_1$, defined in this way are continuous and, due to the formulae (19), both of them are uniformly contracting.

*4.9. Symbolic dynamics*

Consider the set $\Sigma$, consisting of infinite one-side sequences $a = \{a_k \in \{0,1\} : k \in \mathbb{Z}^+\}$. Let us define the metric in the set by a standard way

$$d(a,b) = \sum_{k=0}^\infty 2^{-k}|a_k - b_k|. \quad (21)$$

We identify the periodic subsequences of $\Sigma$, which can be obtained by infinite repetition of finite sequences

$$(a_0, \ldots, a_N),$$

$(a_j \in \{0,1\})$ with these finite sequences.

For any periodic sequence $a$, of this type, generated by the set $\{a_0, \ldots, a_N\}$ we put into correspondence the admissible disk $Q_a$, which is the unique fixed point of the mapping

$$G_{a_0} \circ G_{a_1} \circ \ldots \circ G_{a_N}.$$

Then the following statements hold:

$$Q \in \mathcal{Q}_{a_0}, \quad Q \in G_{a_0}(\mathcal{Q}_{a_1}), \quad Q \in G_{a_0} \circ G_{a_1}(\mathcal{Q}_{a_2}), \ldots,$$
$$Q \in G_{a_0} \circ G_{a_1} \circ \ldots \circ G_{a_{N-1}}(\mathcal{Q}_{a_N}). \quad (22)$$

This means that there exists a constant $C > 0$ such that for any $k \in \mathbb{N}$ if the first $k$ elements of finite sequences $a = (a_0, a_1, \ldots, a_{N_1})$ and $b = (b_0, b_1, \ldots, b_{N_2})$ coincide then

$$\text{dist}_1(Q_a, Q_b) \leq C 2^{-k}.$$

Consequently, for any converging sequence $\{b^k\} \subset \Sigma$, consisting of periodic elements, the corresponding sequence $Q_{b^k}$ converges at $\mathcal{Q}$. Fix an element

$$a = (a_0, a_1, \ldots, a_N, \ldots) \in \Sigma$$

and denote $a^k = (a_0, \ldots, a_k)$. Clearly, $a^k \to a$. Let

$$Q_a = \lim Q_{a^k}.$$

Consider

$$K = \{Q_a : a \in \Sigma\}.$$

Let $H$ be a parameterizing mapping: $H(a) = Q_a$ and $\sigma_i$ be the adding of $i \in \{0, 1\}$ to the left side of the element $a \in \Sigma$. For an arbitrary $a \in \Sigma$ we denote $1a = \sigma_1(a)$, $0a = \sigma_0(a)$.

Let us prove that the disks, corresponding to different elements of the space $\Sigma$, are different. Let $a, b \in \Sigma$ be such that, $a \neq b$. Consider the least integer $j$, such that $a_j \neq b_j$. If $j = 0$, the disks $Q_a$ and $Q_b$ appertain to different sets $U_i$. Otherwise, due to (22), we have

$$Q_a \in G_{a_0} \circ G_{a_1} \circ \ldots \circ G_{a_{j-1}}(\mathcal{Q}_{a_j}),$$
$$Q_b \in G_{b_0} \circ G_{b_1} \circ \ldots \circ G_{b_{j-1}}(\mathcal{Q}_{b_j}) =$$
$$G_{a_0} \circ G_{a_1} \circ \ldots \circ G_{a_{j-1}}(\mathcal{Q}_{b_j}).$$

Hence these disks do not intersect.

**Lemma 4.** *For both $i \in \{0, 1\}$*

$$\sigma_i \circ H = H \circ G_i. \tag{23}$$

*or, equivalently, $D_{ia} = G_i(D_a)$ for any $a \in \Sigma$, $i \in \{0, 1\}$.*

**Proof.** Let $i = 0$, the case $i = 1$ is similar. Fix a sequence $a = \{a_k : k \in \mathbb{Z}^+\} \in \Sigma$ and the corresponding disk

$$Q_a = \mathcal{Q}_{a_0} \bigcap G_{a_0}(\mathcal{Q}_{a_1}) \bigcap G_{a_0} \circ G_{a_1}(\mathcal{Q}_{a_2}) \bigcap \ldots$$
$$\bigcap G_{a_0} \circ G_{a_1} \circ \ldots \circ G_{a_{N-1}}(\mathcal{Q}_{a_N}) \bigcap \ldots$$

Then

$$Q_{0a} = \mathcal{Q}_0 \bigcap G_0(\mathcal{Q}_{a_0}) \bigcap G_0 \circ G_{a_0}(\mathcal{Q}_{a_1}) \bigcap \ldots$$
$$\bigcap G_0 \circ G_{a_0} \circ G_{a_1} \circ \ldots \circ G_{a_{N-1}}(\mathcal{Q}_{a_N}) \bigcap \ldots = G_0(Q_a).$$

∎

The periodic points of the shift mapping are dense in $\Sigma$, and this set is transitive i.e. there is a point

$$a^* = \{a_k^*, k \in \mathbb{Z}^+\} \in \Sigma,$$

whose shifts are dense [35]. Consequently, the periodic admissible disks are dense in $K$ and there is a dense sequence $\{Q_k : k \in \mathbb{Z}^+\} \in K$ such that

$$G_{a_k^*}(Q_k) = Q_{k-1}$$

for all $k \in \mathbb{N}$. Here all the disks $Q_k$ correspond to shifts of the sequence $a_k^*$.

Note, that the similar set of admissible disks exists for any mapping $J$, satisfying (20).

### 4.10. The main result

Recall some notations: $S_{\mu,\theta}$ is the Poincaré period shift mapping for the considered vibro-impact system (2), $z_{\mu,\theta}$ is the fixed point of this mapping, corresponding to the grazing family of periodic solutions,

$$A = \lim_{\theta \to 0+} \lim_{\mu \to 0+} \frac{\partial z}{\partial z_0}(T + \theta - 0, \theta + 0, z_0, \mu)|_{z_0 = z_{\mu,\theta}}.$$

We have proved the following statement.

**Theorem 2.** *Suppose the vibro-impact system (2) satisfies the Condition 2 (existence of the grazing family) and the elements of the matrix $A$ satisfy (16). Then there exists numbers $\mu^* > 0$, $\theta^* > 0$ such that for any $\mu \in (0, \mu^*)$, $\theta \in (0, \theta^*)$ the domains $U_i$, defined by (18) and the corresponding domains $V_i$, $i = 0, 1$ are such that the following statements are true.*

1. *There is an integer $k$ and domains $V_0 \subset U_0$, $V_1 \subset U_1$, such that for any disk $D$, admissible at $U_0$ or $U_1$, (the definition of admissible disk is given at the first paragraph of the Subsection 4.7), the set $S_{\mu,\theta}^k(D)$ contains a disk, admissible at $V_0$ and a disk, admissible at $V_1$.*

2. *The set $\mathcal{J}$ of Lipschitz (with the constant, equal to 1) embeddings $J : \mathcal{Q}^- \to \mathcal{Q}$ is not empty. Here $\mathcal{Q}^-$ and $\mathcal{Q}$ are the sets of disks, admissible at $V_i$ and $U_i$ respectively.*

3. *For any $J \in \mathcal{J}$ there exists a continuous embedding $H : \Sigma \to \mathcal{Q}$, of the space $\Sigma$ consisting of one side boolean sequences with the metric (21) to the space $\mathcal{Q}$ with the $C^1$ metric. This embedding conjugates mappings $\sigma_i$ (adding of $0$ or $1$ to a sequence) and mappings $G_i$ in the sense of formulae (23).*

**Remark 8.** A non-formal reformulation of this statement is given at the Section 6.

### 4.11. Structural stability

Let us discuss the robustness of the constructed set with respect to small perturbations of the diffeomorphism

$S = S_{\mu,\theta}$ or, equivalently, the mapping $\hat{S} = \hat{S}_{\mu,\theta}$. Fix a number $\delta > 0$ and consider a diffeomorphism $\widetilde{S}$ such that

$$|\hat{S}(x) - \widetilde{S}(x)| \leq \delta, \qquad |D\hat{S}(x) - D\widetilde{S}(x)| \leq \delta$$

for any $x \in U_0 \bigcup U_1$. We suppose the number $\delta$ be so small that for any $Q \in \mathcal{Q}$ the image $\widetilde{S}^M(Q)$ contains two admissible subdisks of the sets $V_i$. Also, we suppose that

$$|\Pi_1 \hat{S}_{\mu,\theta}(\zeta_1^s, \zeta_{tan}) - \Pi_1 \widetilde{S}(\zeta_2^s, \zeta_{tan})| \leq \frac{2}{3} |\zeta_1^s - \zeta_2^s|$$

for all values $(\zeta_i, \zeta_{tan}) \in U$.

Then, similarly to what we have done above, we may construct the set $\widetilde{K}$ of admissible disks, with the properties, similar to ones of the set $K$.

## 5. Example

Consider a two degree-of-freedom system, given by equations

$$\ddot{x}_1 + 2p\dot{x}_1 + qx_1 - x_2 = 0; \qquad \ddot{x}_2 + \omega^2 x_2 = a. \tag{24}$$

Suppose that

$$q - p^2 = \omega_0^2 > 0, \quad a > 0, \quad \omega > 0, \quad p^2 + (q-\omega)^2 > 0 \tag{25}$$

and there exists $k \in \mathbb{N}$, such that

$$\omega_0/\omega \in (k + 1/4, k + 1/2). \tag{26}$$

Denote $T = 2\pi/\omega$. Suppose that the system (24) is defined for $x_1 \geq 0$. The impact conditions are the following.
**Condition 3.**

1. Let $r \in (0, 1]$. If $x_1(t_0) = 0$, then

$$x_1(t_0 + 0) = 0, \quad x_2(t_0 + 0) = x_2(t_0 - 0),$$

$$\dot{x}_1(t_0 + 0) = -r\dot{x}_1(t_0 - 0), \quad \dot{x}_2(t_0 + 0) = \dot{x}_2(t_0 - 0).$$

2. Let the instant $t_0$ be such that $x_1(t_0) = \dot{x}_1(t_0 - 0) = 0$ and $x_2(t_0) \leq 0$. Consider the solution $\xi(t)$ of the second equation (24) with initial conditions

$$\xi(t_0) = x_2(t_0 - 0), \quad \dot{\xi}(t_0) = \dot{x}_2(t_0 - 0).$$

Let $I = [t_0, t_1]$ be the maximal segment such that $x_2(t) \leq 0$ for all $t \in I$. Then $x_1(t)|_I \equiv 0$.

In this section the following VIS is considered:

$$\begin{gathered} \ddot{x}_1 + 2p\dot{x}_1 + qx_1 - x_2 = 0; \qquad \ddot{x}_2 + \omega^2 x_2 = a, \\ \text{while } x_1(t) > 0; \\ \text{the Condition 3 is applied if } x_1(t) = 0. \end{gathered} \tag{27}$$

Note that

$$x_2(t) = \frac{a}{\omega^2} + b\sin(\omega(t + \varphi_0)).$$

Here $b$ and $\varphi_0$ are constants. Then the system (24) may be reduced to the equation

$$\ddot{x}_1 + 2p\dot{x}_1 + qx_1 = \frac{a}{\omega^2} + b\sin(\omega(t + \varphi_0)). \tag{28}$$

Provided the conditions (25) are satisfied, the Eq. (28) has the periodic solution

$$x_1^*(t) = \\ \frac{a}{\omega^2 q} + \frac{(q - \omega^2)b\sin(\omega(t + \varphi_0))}{4p^2\omega^2 + (q - \omega^2)^2} - \frac{2pb\omega \cos(\omega(t + \varphi_0))}{4p^2\omega^2 + (q - \omega^2)^2}.$$

Let $b^*(a) = \sqrt{4p^2\omega^2 + (q - \omega^2)^2} a/(q\omega^2)$. If $b = b^*(a)$, the function $x_1^*(t)$ is non-negative and there is exactly one impact instant of the corresponding solution over the period $[0, 2\pi/\omega)$, which is grazing. The general solution of the Eq. (28) is of the form

$$x_1(t) = C \exp(-pt) \sin(\omega_0(t + \theta)) + \frac{a}{\omega^2 q} + \\ \frac{(q - \omega^2)b\sin(\omega(t + \varphi_0)) - 2pb\omega \cos(\omega(t + \varphi_0))}{4p^2\omega^2 + (q - \omega^2)^2}.$$

Here $C$ and $\theta$ are constants. All the (1,1) periodic solutions of the system (27) are ones of the boundary problem, depending on the parameter $t_0 \in [0, T)$

$$\begin{cases} x(t_0) = x(t_0 + T) = 0; & \dot{x}(t_0) = -r\dot{x}(t_0 + T); \\ x(t) > 0 \quad \text{for } t \in (t_0, t_0 + T). \end{cases} \tag{29}$$

Fix $b$ and $\varphi_0$. The Eq. (29) may be considered as a system on the parameters $C$, $\theta$ and $t_0$. Denote

$$\begin{gathered} t_1 = \mathbb{R}csin(2p\omega/\sqrt{4p^2\omega^2 + (q - \omega)^2}), \\ \tau_0 = t_0 + \varphi_0 - t_1, \qquad \tau = t + \varphi_0 - t_1, \\ A = a/(\omega^2 q), \quad B = b/\sqrt{4p^2\omega^2 + (q - \omega)^2}), \\ C_1 = C \exp(pt_0), \quad \vartheta = \theta + t_0. \end{gathered}$$

The general solution of the Eq. (28) is of the form

$$x_1(\tau) = C_1 \exp(-p(\tau - \tau_0)) \sin(\omega_0(\tau - \tau_0 + \vartheta)) + \\ A + B\sin(\omega\tau).$$

The first two equations of (29) may be rewritten as follows:

$$\begin{cases} C_1 \sin(\omega_0 \vartheta) + A + B \sin(\omega \tau_0) = 0; \\ C_1 \sin(\omega_0 \vartheta) = C_1 \exp(-pT) \sin(\omega_0(\vartheta + T)); \\ C_1(\omega_0 \cos(\omega_0 \vartheta) - p \sin(\omega_0 \vartheta)) + B\omega \cos(\omega \tau_0) = \\ \quad -rC_1 \exp(-pT)(\omega_0 \cos(\omega_0(\vartheta + T)) - \\ \quad p\sin(\omega_0(\vartheta + T))) - rB\omega \cos(\omega(\tau_0 + T)). \end{cases} \quad (30)$$

The second of the Eq. (30) gives two cases: either $C_1 = 0$, or

$$\cot(\omega_0 \vartheta) = \frac{\exp(pT) - \cos(\omega_0 T)}{\sin(\omega_0 T)}. \quad (31)$$

If $C_1 = 0$, then $b = b^*(a)$ and function $x_1(t)$ has exactly one zero (of the multiplicity 2) over the period. Otherwise the condition (31) defines uniquely the number

$$\vartheta \in (0, \pi/\omega_0).$$

It follows from the third equation (30) that

$$C_1 = B(1+r)\omega \cos(\omega \tau_0)/H. \quad (32)$$

Here

$$H = r\exp(-pT)(p\sin(\omega_0(\vartheta + T)) - \\ \omega_0 \cos(\omega_0(\vartheta + T))) + p\sin(\omega_0 \vartheta) - \omega_0 \cos(\omega_0 \vartheta).$$

Substituting the expression (32) to the first equation (30), we obtain

$$B(\sin(\omega \tau_0) + D\cos(\omega \tau_0)) + A = 0, \quad (33)$$

where $D = (1+r)\omega \sin(\omega_0 \vartheta)/H$. The number $\tau_0$ may be found from the Eq. (33) if and only if

$$B^2(1 + D^2) \geq A^2. \quad (34)$$

This inequality implies $|b| \geq b^*(a)$. If (34) is false, the system (29) is unsolvable. For $b = \pm b^*(a)$ the graph of the function $x_1(t)$ is tangent to $Ox$ and does not have any other zeros over the period. If $|b| > b^*(a)$, then for any fixed value of the pair $(b, \varphi_0)$ there are two solutions of the Eq. (33) and, respectively, two periodic solutions of the VIS (27). These solutions continuously depend on $b$ and $\varphi_0$ and have a single impact over the period. The corresponding velocities tend to 0 as $b \to \pm b^*(a)$. The system

$$\begin{cases} \dot{x} = y, \\ \dot{y} = -qx - 2py \end{cases} \quad (35)$$

is the homogeneous part for (28). Let $\Phi(t)$ be the fundamental matrix of the system (35), turning to the unity at $t = 0$. The trace of $\Phi(T)$ equals to $2\cos(\omega_0 T) \exp(-pT)$, and the determinant equals to $\exp(-2pT)$. Due to the inequality (26), the conditions of the Theorem 2 are satisfied and there exists a set of admissible disks, corresponding to the VIS (27), which can be described by means of symbolic dynamics.

## 6. Discussion.

Here we reformulate the result of the Theorem 2 on a less formal level. We try to find a "visible" analogue of the chaotic dynamics in a neighborhood of an invisible chaotic invariant set. First of of all we construct a non-hyperbolic homoclinic point.

For this we need conditions (16). These conditions have a clear geometric sense. The inequality $a_{12} > 0$ implies that the unstable manifold $W^u$ exists and intersects the grazing surface $\Gamma_{\mu,\theta}$ provided $\mu$ and $\theta$ are positive and sufficiently small. The unstable manifold bents at the point of the intersection. If

$$\sum_{j=1}^{2n} a_{1k} a_{k2} < -a_{12}$$

the manifold $W^u$ turns so that it intersects $W^{cs}$ transversally.

This does not imply the existence of a Smale horseshoe. What we can expect is to find a symbolic dynamics on leaves of a foliation over a Cantor set. However, even this is not always possible.

A central stable manifold is not as "good" as an ordinary stable manifold. For example, it may be non-unique even for a fixed point [44]. Typically, a central stable bundle over an invariant set is not integrable [43] i.e. we can not find a corresponding continuous foliation, like what we can do in the hyperbolic case. Moreover, the iterations of a central stable (or a central unstable) manifold may shrink both in positive and negative direction. So, we have to be able to expand iterated leaf on every step. This is why we need the mapping $J$ from the statement of the theorem.

If we can reconstruct the whole central stable (or central unstable) leave by any small plaque of this disk (this is the so-called *dynamical coherence* [45]), then we get a naturally defined expanding mapping $J$ and, due to the Theorem 2, a symbolic dynamics on the set of central unstable leaves. For example, this happens, if the considered central unstable manifolds are de facto ordinary unstable manifolds (a hyperbolic case) or, more generally, there is a Lyapunov function, which does not allow the central unstable leaves to shrink. For the last case, we do not need any expanding mappings at all.

In the general case, we offer a sample of expanding mapping (Subsection 4.8) drawing straight lines through the boundary points of the "narrow" disk. The expansion procedure allows us to construct a foliation. So, iterating and expanding disks on every step, we obtain a complex dynamics on the set of disks, close to the central unstable bundle. We have constructed a model of a slow-fast system with a chaotic dynamics for the fast variable and unknown behavior of the slow one.

Compare the Devaney's definition of chaos with one offered in our paper.

**Remark 9.**

1. The Devaney's definition of chaos requests hyperbolicity. The introduced model does not.

2. The Devaney's definition describes the dynamic of point. Our model describes the dynamic of approximations to central unstable leaves.

3. Both the considered models are structurally stable.

## 7. Conclusion.

An impact oscillator, described by the Newtonian model of impacts has been studied. There exist many results on grazing bifurcations for such systems (see [10-12] for review). It has been shown that different types of chaos can be observed for near-grazing values of parameters. The explicit conditions, providing existence of a complex dynamics have been given. The geometry of the obtained strange attractors has been described.

However, especially in the case of several degrees-of-freedom, the structure of these attractors may be very sensitive to the parameters of the system. This happens because the ratios of the Lyapunov exponents of this solution are big in a non-degenerate case.

In order to describe the robust properties of attractors, we offer a new approach treating the near-grazing dynamics as non-hyperbolic. We use a technique of central unstable manifolds. A dynamical system on a set of embeddings of the local central unstable manifold of a periodic point is introduced. For this system we find an invariant set described by the symbolical dynamics.

The conditions of existence of this set (Theorem 2) do not require hyperbolicity of the periodic motion. So, the obtained structure is less sensitive to the changes of parameters when the strange attractor, which may disappear if one of the Lyapunov exponents becomes zero. Roughly speaking, some traces of the strange attractor persist, even when the attractor itself has already disappeared.

From the mathematical point of view, we do not use any common assumptions from the theory of partial hyperbolicity, like integrability of the central unstable bundle, etc. There may be no continuous central unstable foliation in the considered case. What we have constructed, is not a classical "dynamics of leaves" [43]. However, all the "admissible disks" are close to the central unstable bundle and can be considered as approximations to local central unstable manifolds. So, the result of the Theorem 2 describes a new type of chaotic dynamics.

Grace to the work by Dankowitz and Nordmark [7], the results of the current paper may be transferred to the "soft" models of impacts.

The obtained model of near-grazing dynamics is structurally stable and can be observed in experiments.

**Acknowledgements.**

This work was supported by the UK Royal Society, by the Russian Federal Program "Scientific and pedagogical cadres", grant no. 2010-1.1-111-128-033 and by the Chebyshev Laboratory (Department of Mathematics and Mechanics, Saint-Petersburg State University) under the grant 11.G34.31.2006 of the Government of the Russian Federation.

## References

[1] F. Peterka *Laws of impact motion of mechanical systems with one degree of freedom* Acta Technika CSAV, 1974, 1−11.

[2] J. M. T. Thompson, R. Ghaffari *Chaotic dynamics of an impact oscillator.* Phys. Rev. A. **27** (1983), 1741−1743.

[3] A. B. Nordmark *Non-periodic motion caused by grazing incidence in an impact oscillator.* J. Sound Vib., **145** (1991), 279–297.

[4] *Nordmark A. B.* Effects due to low velocity impact in mechanical oscillators // Int. J. Bifur. Chaos Appl. Sci. Engrg. V. 2, 1992, P. 597−605.

[5] *Nordmark A. B.* Universal limit mapping in grazing bifurcations // Phys. Rev. E. V. 55, 1997. P. 266−270.

[6] M. H. Fredriksson, A. B. Nordmark *Bifurcations caused by grazing incidence in many degrees of freedom impact oscillators.* Proc. Roy. Soc. London Ser. A. **453** (1997), 1261–1276.

[7] H. Dankowicz, A. B. Nordmark *On the origin and bifurcations of stick-slip oscillators*, Physica D, **136** (1999), 280–302.

[8] A. B. Nordmark *Existence of periodic orbits in grazing bifurcations of impacting mechanical oscillators.* Nonlinearity **14** (2001), 1517−42.

[9] H. Dankowicz, P. Piiroinen, A. B. Nordmark *Low-velocity impacts of quasiperiodic oscillations.* Chaos Solitons Fractals **14** (2002), 241−255.

[10] *di Bernardo M., Budd C. J., Champneys A. R., Kowalczyk P., Nordmark A. B., Olivar G., Piirionen P. T.* Bifurcations in Nonsmooth Dynamical Systems SIAM Review. **50**, 2008. P. 629−701.

[11] M. di Bernardo, C. J. Budd, A. R. Champneys, P. Kowalczyk *Bifurcations and Chaos in Piecewise-Smooth Dynamical Systems: Theory and Applications.* New York, Springer, 2007.

[12] M. Kunze *Non-Smooth Dynamical Systems* Lecture Notes in Mathematics, **1744**, 2000, 228 p.

[13] G. S. Whiston *The vibro-impact response of a harmonically excited and preloaded one-dimensional linear oscillator* J. Sound Vib. **115** (1987), 303−319.

[14] G. S. Whiston *Global dynamics of a vibro-impacting linear oscillator* J. Sound Vib. **118** (1987), 395−429.

[15] G. S. Whiston G. S. *Singularities in vibro-impact dynamics* J. Sound Vib., **152** (1992), 427−460.

[16] C. J. Budd, F. Dux *Chattering and related behavior in impacting oscillators* Phil. Trans. Roy. Soc. **347** (1994), 365−389.

[17] C. J. Budd, F. Dux *Intermittency in impact oscillators close to resonance* Nonlinearity, **7** (1994), 1191−1224.

[18] D. R. J. Chillingworth *Discontinuity geometry for an impact oscillator* Dynamical Systems, **17** (2002), 389–420.

[19] D. R. J. Chillingworth *Dynamics of an impact oscillator near a degenerate graze* Nonlinearity, **23** (2010), 2723-2748.

[20] H. M. Osinga, R. Szalai *Unstable manifolds of a limit cycle near grazing.* Nonlinearity **21** (2008) 273−284.

[21] R. L. Devaney *An Introduction to Chaotic Dynamical Systems.* Redwood City, CA: Addison-Wesley, 1987.

[22] J. Molenaar, W. van de Water, J. de Wegerand *Grazing impact oscillations.* Phys. Rev. E, **62** (2000), 2030–2041.

[23] E. Pavlovskaia, M. Wiercigroch *Low-dimensional maps for piecewise smooth oscillators* J. Sound Vib. **305**, 2007, 750−771.

[24] S. Banerjee, J. Ing J., E. Pavlovskaia, M. Wiercigroch *Experimental study of impact oscillator with one-sided elastic constraint* Phil. Trans. R. Soc. A., **366** (2008), 679704.

[25] S. Banerjee, Ing J., Pavlovskaia E., Wiercigroch M., *Invisible grazings and dangerous bifurcations in impacting systems: The problem of narrow-band chaos* Phys. Rev. E **79** (2009) 037201.

[26] J. Ing, E. Pavlovskaia, M. Wiercigroch, S. Banerjee *Bifurcation analysis of an impact oscillator with a one-sided elastic constraint near grazing.* Physica D 239 (2010) 312–321.

[27] H. E. Nusse, E. Ott, J. A. Yorke *Border-collision bifurcations: An explanation for observed bifurcation phenomena* Phys. Rev. E., **49** (1994), 1073−1076.

[28] S. Banerjee, J. A. Yorke, C. Grebogi *Robust chaos.* Phys. Rew. Letters. **80**, (1998), no. 14, 3049–3052.

[29] R. Szalai, G. Stepan, S. J. Hogan, *Global dynamics of low immersion high-speed milling.* Chaos, **14** (2004), Issue 4, 1069–1077.

[30] W. Chin, E. Ott, H. E. Nusse, C. Grebogi *Universal behavior of impact oscillators near grazing incidence.* Phys. Let. A. **201**, (1995), 197–204.

[31] C. Gontier, C. Toulemonde *Approach to the periodic and chaotic behaviour of the impact oscillator by a continuation method* European J. Mech. A. Solids. **16** (1997), 141−163.

[32] P. J. Holmes *The dynamics of repeated impacts with a sinusoidally vibrating table.* J. Sound. Vib. **84**, (1982), 173–189.

[33] S. W. Shaw, R. H. Rand *The transition to chaos in a simple mechanical system.* Int. J. Non-Linear Mechanics, **24** (1989), 41−56.

[34] S. G. Kryzhevich *Grazing bifurcation and chaotic oscillations of single degree-of-freedom dynamical systems.* J. Appl. Math. Mech. **72** (2008), 539–556.

[35] S. Smale, *Diffeomorphisms with many periodic points*, Differential and Combinatorial Topology, S. S. Caimes, ed. (Princeton University Press, Princeton, New Jersey, 1965), pp. 68-80.

[36] M. U. Akhmet *Li-Yorke chaos in systems with impacts.* J. Math. Anal. Appl. **351**, (2009), 804–810.

[37] S. P. Gorbikov, A. V. Men'shenina *Bifurcation, leading to chaotic motions in dynamical systems with impact interactions.* Diff. equations, **41** (2005), 1046−1052.

[38] S. P. Gorbikov, A. V. Men'shenina *Statistical description of the limiting set for chaotic motion of the vibro-impact system.*, Automation and remote control, **68** (2007), 1794–1800.

[39] Bo Deng *The Shil'nikov problem, exponential expansion, strong $\lambda$-lemma, $C^1$-linearization, and homoclinic bifurcation.* J. Diff. Equations 79 Issue 2 (1989), 189–231.

[40] G. Cicogna, M. Santoprete *Nonhyperbolic homoclinic chaos* Physics Letters A 256 (1999), 25–30.

[41] A. Gorodetski, Yu. S. Ilyashenko *Some properties of skew products over the horseshoe and solenoid* Proceedings of Steklov Institute 231 (2000), 96–118.

[42] S. V. Gonchenko, L. P. Shil'nikov, O. V. Sten'kin, D. V. Turaev *Bifurcations of Systems with Structurally Unstable Homoclinic Orbits and Moduli $\Omega$-Equivalence* Computer Mathematics and Applications. 34 (1997), no. 24, 111–142.

[43] M. W. Hirsch, C. C. Pugh, M. Shub *Invariant Manifolds*, Springer-Verlag, Berlin-Heidelberg, 1977. 154 pp.

[44] A. Kelley *The Stable, Center Stable, Center, Center Unstable, Unstable Manifolds* J. Diff. Equations **3** (1967), 546−570.

[45] Ya. B. Pesin *Lectures on partial hyperbolicity and stable ergodicity*, Zurich Lectures in Advanced Mathematics (9783037190036) 2006.

[46] V. A. Pliss *Reduction principle in the stability theory* Izv. AN SSSR, Mat. ser. **28** (1964), 1297−1324.

[47] A. Reinfelds *A reduction principle for discrete dynamical and semidynamical systems in metric spaces* Z. Agnew Math. Phys.

[48] A. P. Ivanov *Stabilization of an impact oscillator near grazing incidence owing to resonance* J. Sound Vib. **162**, (1993), no. 3, 562−565.

[49] A. P. Ivanov *Bifurcations in impact systems.* Chaos Solitons Fractals **7** (1996), 1615–34. **45** (1994), 933−955.